\RequirePackage{lineno}
\documentclass[12pt,a4paper,leqno,verbatim]{amsart}

\newcounter{minutes}
\divide\time by 60
\newcounter{hours}
\multiply\time by 60 \addtocounter{minutes}{-\time}

\usepackage{amssymb}
\usepackage{hyperref}
\usepackage[T1]{fontenc}
\usepackage{graphicx}
\date{}
\newfont{\cyrilic}{wncyr10 scaled 1000}

\title[generalized trigonometric and hyperbolic functions]
{Inequalities for the generalized trigonometric and hyperbolic functions}

\author{Riku Kl\'en}
\author{Matti Vuorinen}
\author{Xiaohui Zhang }

\address{Department of Mathematics and statistics, University of Turku, 20014 Turku,
Finland} \email{ripekl@utu.fi, vuorinen@utu.fi, xiazha@utu.fi}

\newcommand{\comment}[1]{}

\swapnumbers
\theoremstyle{plain}

\newtheorem{theorem}[equation]{Theorem}
\newtheorem{lemma}[equation]{Lemma}

\newtheorem{corollary}[equation]{Corollary}

\newtheorem{conjecture}[equation]{Conjecture}
\newtheorem{openproblem}[equation]{Open problem}

\newcommand{\R}{\mathbb{R}}
\newcommand{\Z}{\mathbb{Z}}

\newcommand{\snp}{\sin_p}
\newcommand{\cp}{\cos_p}
\newcommand{\tp}{\tan_p}
\newcommand{\asp}{\arcsin_p}

\newcommand{\atp}{\arctan_p}
\newcommand{\shp}{\sinh_p}
\newcommand{\chp}{\cosh_p}
\newcommand{\thp}{\tanh_p}
\newcommand{\ashp}{{\rm{arcsinh}}_p}

\newcommand{\athp}{{\rm{arctanh}}_p}
\newcommand{\pip}{\pi_p}

\newcommand{\beq}{\begin{equation}}
\newcommand{\eeq}{\end{equation}}
\newcommand{\bthm}{\begin{theorem}}
\newcommand{\ethm}{\end{theorem}}
\newcommand{\bp}{\begin{proof}}
\newcommand{\ep}{\end{proof}}
\newcommand{\blem}{\begin{lemma}}
\newcommand{\elem}{\end{lemma}}
\newcommand{\bcol}{\begin{corollary}}
\newcommand{\ecol}{\end{corollary}}

\numberwithin{equation}{section}

\begin{document}

\begin{abstract}
The generalized trigonometric functions occur as an eigenfunction of the Dirichlet problem for the one-dimensional $p-$Laplacian. The generalized hyperbolic functions are defined similarly. Some classical inequalities for trigonometric and hyperbolic functions, such as Mitrinovi\'c-Adamovi\'c inequality, Lazarevi\'c's inequality, Huygens-type inequalities, Wilker-type inequalities, and Cuza-Huygens-type inequalities, are generalized to the case of generalized functions.
\end{abstract}

\def\thefootnote{}
\footnotetext{ \texttt{\tiny File:~\jobname .tex,
          printed: \number\year-\number\month-\number\day,
          \thehours.\ifnum\theminutes<10{0}\fi\theminutes}
} \makeatletter\def\thefootnote{\@arabic\c@footnote}\makeatother

\maketitle




\section{Introduction}

It is well known from basic calculus that
$$
{\rm arcsin}(x)=\int_0^x \dfrac{1}{(1-t^2)^{1/2}}dt, \quad 0\leq x\leq1,
$$
and
$$
\dfrac{\pi}{2}={\rm arcsin}(1)=\int_0^1 \dfrac{1}{(1-t^2)^{1/2}}dt.
$$
We can define the function ${\rm sin}$ on $[0,\pi/2]$ as the inverse of ${\rm arcsin}$ and extend it on $(-\infty,\infty)$.

Let $1<p<\infty$. We can generalize the above functions as follows:
$$
\asp(x)\equiv\int_0^x \dfrac{1}{(1-t^p)^{1/p}}dt, \quad 0\leq x\leq1,
$$
and
$$
\dfrac{\pip}{2}=\asp(1)\equiv\int_0^1 \dfrac{1}{(1-t^p)^{1/p}}dt.
$$
The inverse of $\asp$ on $[0,\pip/2]$ is called the \emph{generalized sine function} and denoted by $\snp$.
By standard extension procedures as the sine function we get a differentiable function on the whole of
$(-\infty,\infty)$ which coincides with ${\rm sin}$ when $p=2$.
It is easy to see that the function $\snp$ is strictly increasing and concave on $[0,\pip/2]$.
In the same way we can define the generalized cosine function, the generalized tangent function, and their inverses, and also the corresponding hyperbolic functions.

The generalized sine function $\snp$ occurs as an eigenfunction of the Dirichlet problem for the one-dimensional $p-$Laplacian. There are several different definitions for these generalized trigonometric and hyperbolic functions
\cite{le, l1, l2, lp1}.
Recently, these functions have been studied very extensively (see \cite{bv, be, egl, le, l1, l2, lp1, lp2}). In particular, the reader is referred to \cite{l1, l2, lp1, lp2}. These generalized functions are similar to the classical functions in various aspects. Some of these functions can be expressed in terms of the Gaussian hypergeometric series (see \cite{bv1, bv}).

In this paper we will generalize some classical inequalities for trigonometric and hyperbolic functions, such as \emph{Mitrinovi\'c-Adamovi\'c inequality} (Theorem \ref{mit-ada}), \emph{Lazarevi\'c's inequality} (Theorem \ref{laza-p}), \emph{Huygens-type inequalities} (Theorem \ref{thm:huygen1} and Theorem \ref{thm:huygen2} ), \emph{Wilker-type inequality} (Corollary
\ref{thm:wilker}), and \emph{Cuza-Huygens-type inequalities} (Theorem \ref{cusa-hyugens4trig} and Theorem \ref{cusa-hyugens4hyp}) to the case of generalized functions. For the classical cases, these inequalities have been extended and sharpened extensively (see the very recent survey \cite{avz}).

\section{Definitions and formulas}

In this section we define the generalized cosine function, the generalized tangent function, and their inverses, and also the corresponding hyperbolic functions.

The \emph{generalized cosine function} $\cp$ is defined as
$$
\cp(x)\equiv\frac{d}{dx}\snp(x).
$$
It is clear from the definitions that
$$
\cp(x)=(1-\snp(x)^p)^{1/p},\qquad x\in[0,\pip/2],
$$
and
\beq\label{idty:c-s}
|\snp(x)|^p+|\cp(x)|^p=1, \qquad x\in\R.
\eeq
It is easy to see that
$$
\frac{d}{dx}\cp(x)=-\cp(x)^{2-p}\snp(x)^{p-1},\qquad x\in[0,\pip/2].
$$

The \emph{generalized tangent function} is defined as in the classical case:
$$
\tp(x)\equiv\dfrac{\snp(x)}{\cp(x)},\qquad x\in\R\setminus\{k\pip+\frac{\pip}{2}:\,k\in \Z\}.
$$
It follows from (\ref{idty:c-s}) that
$$
\frac{d}{dx}\tp(x)=1+|\tp(x)|^p,\quad x\in(-\pip/2,\pip/2).
$$

Similarly, the \emph{generalized inverse hyperbolic sine function}
$$
\ashp(x)\equiv
\left\{
\begin{array}{ll}
\int_0^x\frac{1}{(1+t^p)^{1/p}}dt,& x\in[0,\infty),\vspace{2mm}\\
-\ashp(-x), & x\in(-\infty,0)
\end{array}\right.
$$
generalizes the classical inverse hyperbolic sine function.
The inverse of $\ashp$ is called the \emph{generalized hyperbolic sine function} and denoted by $\shp$.
The \emph{generalized hyperbolic cosine function} is defined as
$$
\chp(x)\equiv\frac{d}{dx}\shp(x).
$$
The definitions show that
$$
\chp(x)^p-|\shp(x)|^p=1, \qquad x\in \R,
$$
and
$$
\frac{d}{dx}\chp(x)=\chp(x)^{2-p}\shp(x)^{p-1},\qquad x\geq0.
$$
The \emph{generalized  hyperbolic tangent function} is defined as
$$
\thp(x)\equiv\frac{\shp(x)}{\chp(x)},
$$
and hence we have
$$
\frac{d}{dx}\thp(x)=1-|\thp(x)|^p.
$$

It is clear that all these generalized functions coincide with the classical ones when $p=2$.

\section{Inequalities}

The l'H\^opital Monotone Rule (LMR), Lemma \ref{thm:LMR}, is the key tool in proofs of our generalizations.

\blem\label{thm:LMR}\cite{avv}
\rm{(l'H\^opital Monotone Rule}).  Let  $-\infty < a < b < \infty$, and let  $f,g: [a,b] \to \mathbb{R}$ be continuous functions that are differentiable on  $(a,b)$, with $f(a)=g(a)=0$ or $f(b)=g(b)=0.$  Assume that  $g'(x)\ne 0$ for each  $x\in (a,b).$  If  $f'/g'$  is increasing (decreasing) on  $(a,b)$, then so is $f/g$.
\elem

Some other applications of the l'H\^opital Monotone Rule (LMR) in special functions one is referred to the survey \cite{avz}.

\blem\label{lem:tan-tanh}
For $p>2$, the function $f(x)\equiv\tp(x)^{p-2}-\thp(x)^{p-2}$ is strictly increasing in $(0,\pi_p/2)$.
\elem

\bp
By differentiation, we have
$$f'(x)=(p-2)(\tp(x)^{p-3}(1+\tp(x)^p)-\thp(x)^{p-3}(1-\thp(x)^p)).$$
For $p\geq3$,
$$f'(x)\geq(p-2)(\tp(x)^{p-3}-\thp(x)^{p-3})>0,$$
since $\tp(x)>\thp(x).$

By the identities $\snp(x)^p+\cp(x)^p=1$ and $\chp(x)^p-\shp(x)^p=1$,
\begin{eqnarray*}
f'(x)&=&(p-2)\left(\dfrac{\snp(x)^{p-3}}{\cp(x)^{2p-3}}-\dfrac{\shp(x)^{p-3}}{\chp(x)^{2p-3}}\right)\\
&\geq&(p-2)\shp(x)^{p-3}\left(\dfrac{1}{\cp(x)^{2p-3}}-\dfrac{1}{\chp(x)^{2p-3}}\right)>0
\end{eqnarray*}
for $p\in[2,3)$ since $\snp(x)<\shp(x)$. This completes the proof.
\ep

\blem\label{lem:c-ch}
For $p>1$, the function $f(x)\equiv\cp(x)\chp(x)$ is strictly decreasing from $(0,\pi_p/2)$ onto $(0,1)$. In particular,
for all $p\in(1,\infty)$ and $x\in(0,\pi_p/2)$,
$$\cp(x)<\dfrac{1}{\chp(x)}.$$
\elem

\bp
After simple computations we get
$$f'(x)=\cp(x)\chp(x)(\thp(x)^{p-1}-\tp(x)^{p-1})<0,$$
which implies that $f$ is strictly decreasing, and hence $\cp(x)\chp(x)<1.$
\ep

\bthm
For $p\in[2,\infty)$ and $x\in(0,\pi_p/2)$,
\beq\label{2-4.4-r}
\dfrac{\snp(x)}{x}<\dfrac{x}{\shp(x)}.
\eeq
\ethm

\bp
Let $f_1(x)\equiv\snp(x)\shp(x)$, $f_2(x)\equiv x^2$ and $f_1(0)=f_2(0)=0$. By simple computations, we have
$$\dfrac{f_1''(x)}{f_2''(x)}=\cp(x)\chp(x)-\dfrac12\snp(x)\shp(x)(\tp(x)^{p-2}-\thp(x)^{p-2})$$
which is strictly decreasing for any $p\geq2$ by Lemma \ref{lem:tan-tanh} and Lemma \ref{lem:c-ch}. Hence the monotonicity of $f_1(x)/f_2(x)$ follows from the
l'H\^opital Monotone Rule, and this implies
$$\dfrac{\snp(x)\shp(x)}{x^2}<1.$$
\ep

The next two theorems generalize the \emph{Mitrinovi\'c-Adamovi\'c inequality} and \emph{Lazarevi\'c's inequality} (see \cite{mit}). For the classical case of  Theorem \ref{laza-p} also see \cite{lwc}.

\bthm\label{mit-ada}
For $p\in(1,\infty)$, the function
$$f(x)\equiv\dfrac{\log(\snp(x)/x)}{\log\cp(x)}$$
is strictly decreasing from $(0,\pip/2)$ onto $(0,1/(1+p))$. In particular,
for all $p\in(1,\infty)$ and $x\in(0,\pip/2)$,
\beq\label{mitrinovic-adamovic}
\cp(x)^\alpha<\dfrac{\snp(x)}{x}<1
\eeq
with the best constant $\alpha=1/(1+p)$.
\ethm

\bp
Write $f_1(x)\equiv\log(\snp(x)/x)$ and $f_2(x)\equiv\log\cp(x)$. Then $f_1(0)=f_2(0)=0$ and,
by simple computations,
$$\dfrac{f_1'(x)}{f_2'(x)}=\dfrac{\tp(x)-x}{x\tp(x)^p}=\dfrac{f_{11}(x)}{f_{22}(x)},$$
with $f_{11}(x)\equiv\tp(x)-x$, $f_{22}(x)\equiv x\tp(x)^p$, and $f_{11}(0)=f_{22}(0)=0.$
$$\dfrac{f_{11}'(x)}{f_{22}'(x)}=\dfrac{1}{1+p\,g(x)}$$
with
$$g(x)\equiv\dfrac{x}{\snp(x)}\dfrac{1}{\cp(x)^{p-1}}$$
which is strictly increasing. By the l'H\^opital  Monotone Rule we see that
$f(x)$ is strictly decreasing. The limiting values follow from l'H\^opital's Rule easily.
\ep

\bthm\label{laza-p}
For $p\in(1,\infty)$, the function
$$f(x)\equiv\dfrac{\log(\shp(x)/x)}{\log\chp(x)}$$
is strictly increasing from $(0,\infty)$ onto $(1/(1+p),1)$. In particular,
for all $p\in(1,\infty)$ and $x\in(0,\infty)$,
\beq\label{lazarevic}
\chp(x)^\alpha<\dfrac{\shp(x)}{x}<\chp(x)^\beta
\eeq
with the best constants $\alpha=1/(1+p)$ and $\beta=1$.
\ethm

\bp
Write $f_1(x)\equiv\log(\shp(x)/x)$ and $f_2(x)\equiv\log\chp(x)$. Then $f_1(0)=f_2(0)=0$ and,
by simple computations,
$$\dfrac{f_1'(x)}{f_2'(x)}=\dfrac{x-\thp(x)}{x\thp(x)^p}=\dfrac{f_{11}(x)}{f_{22}(x)},$$
with $f_{11}(x)\equiv x-\thp(x)$, $f_{22}(x)\equiv x\thp(x)^p$, and $f_{11}(0)=f_{22}(0)=0.$
$$\dfrac{f_{11}'(x)}{f_{22}'(x)}=\dfrac{1}{1+p\,g(x)}$$
with
$$g(x)\equiv\dfrac{x}{\shp(x)}\dfrac{1}{\chp(x)^{p-1}}$$
which is strictly decreasing. By the l'H\^opital Monotone Rule we see that
$f(x)$ is strictly increasing. The limiting values follow from l'H\^opital's Rule easily.
\ep

\bcol
For all $p\in[2,\infty)$ and $x\in(0,\pi_p/2)$,
\beq\label{2-4.4}
\left(\dfrac{x}{\shp(x)}\right)^{1+p}<\dfrac{1}{\chp(x)}<\dfrac{\thp(x)}{x}<\dfrac{\snp(x)}{x}<\dfrac{x}{\shp(x)}.
\eeq
\ecol

\bp
The first inequality of (\ref{2-4.4}) follows by the left side of (\ref{lazarevic}). The second inequality follows by
$\shp(x)/x>1$, while the third by $\snp(x)>\thp(x)$. The last inequality is the inequality (\ref{2-4.4-r})
\ep

\begin{conjecture}
For $p\in[2,\infty)$, the function
$$
f(x)\equiv\frac{\log (x/\snp(x))}{\log (\shp(x)/x)}
$$
is strictly increasing in $(0,\pi_p/2)$.
\end{conjecture}

Next two theorems show the \emph{Huygens-type inequalities} for the generalized trigonometric and hyperbolic functions.

\bthm\label{thm:huygen1}
Let $p>1$. Then the following inequalities hold
\beq\label{huygens:sinp-cosp}
(p+1)\frac{\snp(x)}{x}+\frac{1}{\cp(x)}>p+2 \quad for \quad x\in(0,\pip/2),
\eeq
and
\beq\label{huygens:shp-chp}
(p+1)\frac{\shp(x)}{x}+\frac{1}{\chp(x)}>p+2 \quad for \quad x>0.
\eeq
\ethm

\bp
The well-known weighted arithmetic-geometric inequality states that
$$
ta+(1-t)b>a^tb^{1-t},
$$
for $a,b>0$, $a\neq b$, and $0<t<1$. Putting $t=(p+1)/(p+2)$, $a=\snp(x)/x$, and $b=1/\cp(x)$, and combining
the left side of (\ref{mitrinovic-adamovic}), we have
$$
(p+1)\frac{\snp(x)}{x}+\frac{1}{\cp(x)}>
(p+2)\left(\frac{\snp(x)}{x}\right)^{(p+1)/(p+2)}\left(\frac{1}{\cp(x)}\right)^{1/(p+2)}>p+2.
$$

Similarly, the inequality (\ref{huygens:shp-chp}) follows from the left side of (\ref{lazarevic}).
\ep

\bthm\label{thm:huygen2}
For $p>1$, the following inequalities hold
\beq\label{huygens:sin-tan}
\dfrac{p\snp(x)}{x}+\dfrac{\tp(x)}{x}>1+p,\quad 0<x<\dfrac{\pip}{2},
\eeq
and
\beq\label{huygens:sh-th}
\dfrac{p\shp(x)}{x}+\dfrac{\thp(x)}{x}>1+p,\quad x>0.
\eeq
\ethm

\bp
Let $f(x)\equiv p\snp(x)+\tp(x)-(1+p)x$. After some elementary computations, we get
$$
f'(x)=p\cp(x)+\tp(x)^p-p
$$
and
$$
f''(x)=p\tp(x)^{p-1}(1-\cp(x)+\tp(x)^p)>0,
$$
which implies that $f'(x)>0$ and $f$ is strictly increasing.
Hence we have $f(x)>0$, and the inequality (\ref{huygens:sin-tan}) follows.

Similarly, put $g(x)\equiv p\shp(x)+\thp(x)-(1+p)x$. We have
$$
g'(x)=p\chp(x)-\thp(x)^p-p
$$
and
$$
g''(x)=p\thp(x)^{p-1}(\thp(x)^p+\chp(x)-1)>0,
$$
from which we get $g'(x)>0$, implying $g(x)>0$. This finishes the proof.
\ep

\bcol\label{thm:wilker}
For $p>1$ and $x>0$,
\beq\label{wilker:sh-th}
\left(\frac{\shp(x)}{x}\right)^p+\frac{\thp(x)}{x}>2.
\eeq
\ecol

\bp
The well-known Bernoulli inequality states that, for $a>1$ and $t>0$,
\beq\label{bernoulli}
(1+t)^a>1+at.
\eeq
Setting $t=\shp(x)/x-1$ and $a=p$ in (\ref{bernoulli}), and then combining the inequality (\ref{huygens:sh-th}), we have
$$
\left(\frac{\shp(x)}{x}\right)^p>1+p\left(\frac{\shp(x)}{x}-1\right)>2-\frac{\thp(x)}{x},
$$
which implies (\ref{wilker:sh-th}).
\ep

The inequality (\ref{wilker:sh-th}) is the so-called \emph{Wilker's inequality}. The following Theorem \ref{cusa-hyugens4trig} and \ref{cusa-hyugens4hyp} present the famous \emph{Cusa-Huygens-type inequalities} for the generalized trigonometric and hyperbolic functions, respectively.

\bthm\label{cusa-hyugens4trig}
For $p\in(1,2]$, the following inequalities
\beq\label{3-(3)}
\dfrac{\snp(x)}{x}<\dfrac{\cp(x)+p}{1+p}\leq \dfrac{\cp(x)+2}{3}
\eeq
hold for all $x\in(0,\pip/2]$.
\ethm

\bp
Let $f(x)\equiv x\cp(x)+px-(1+p)\snp(x)$. By differentiation, we have
$$f'(x)=-\cp(x)(x\tp(x)^{p-1}+p)+p\equiv-g(x)+p,$$
and
$$g'(x)=\cp(x)\tp(x)^{p-2}\left((p-1)(x-\tp(x))+(p-2)x\tp(x)^p\right)<0,$$
which implies $g(x)<g(0)=p$ and $f'(x)>0$. Hence $f(x)$ is strictly increasing and $f(x)>f(0)=0$
which implies the inequality (\ref{3-(3)}) .

The second inequality in (\ref{3-(3)}) is clear since $\cp(x)\leq1.$
\ep

\bthm\label{cusa-hyugens4hyp}
For all $x>0$,
\beq\label{2-4.6(1)}
\dfrac{\shp(x)}{x}<\dfrac{\chp(x)+p}{1+p},\quad if\quad p\in(1,2],
\eeq
and
\beq\label{2-4.6(2)}
\dfrac{\shp(x)}{x}<\dfrac{\chp(x)+2}{3},\quad if\quad p\in[2,\infty).
\eeq
\ethm

\bp
Let $f(x)\equiv x\chp(x)+px-(1+p)\shp(x)$. By differentiation, we have
$$f'(x)=\chp(x)(x\thp(x)^{p-1}-p)+p$$
and
$$f''(x)=\chp(x)\thp(x)^{p-2}\left((p-1)(x-\thp(x))+(2-p)x\thp(x)^p\right)>0,$$
which implies $f'(x)>0$. Hence $f(x)$ is strictly increasing, and $f(x)>f(0)=0$
which implies the inequality (\ref{2-4.6(1)}) .

For the inequality (\ref{2-4.6(2)}), let $h(x)\equiv x\chp(x)+2x-3\shp(x)$. By differentiation, we get
$$h'(x)=\chp(x)\left(x\thp(x)^{p-1}-2\right)+2$$
and
\begin{eqnarray*}
h''(x)&=&\chp(x)\thp(x)^{p-2}\left(x\thp(x)^p-\thp(x)+(p-1)x(1-\thp(x)^p)\right)\\
&\geq&\chp(x)\thp(x)^{p-2}\left(x\thp(x)^p-\thp(x)+x(1-\thp(x)^p)\right)\\
&=&\chp(x)\thp(x)^{p-2}\left(x-\thp(x)\right)>0,
\end{eqnarray*}
which implies $h'(x)>h'(0)=0$, and hence $h(x)$ is strictly increasing and $h(x)>h(0)=0$. This implies the inequality (\ref{2-4.6(2)}).
\ep

\bthm
For $p\in[2,\infty)$ and $x\in(0,\pi_p/2)$,
\beq\label{2-4.1-l}
\dfrac{\sinh_p(x)}{x}<\dfrac{3}{2+\cos_p(x)}.
\eeq
\ethm

\bp
Let
$$f(x)\equiv3x-2\sinh_p(x)-\sinh_p(x)\cos_p(x).$$
Simple computations give
\begin{eqnarray*}
f'(x)&=&3-2\chp(x)-\chp(x)\cp(x)+\shp(x)\snp(x)^{p-1}\cp(x)^{2-p}\\
&\geq&3-2\chp(x)-\chp(x)\cp(x)+\shp(x)\snp(x)^{p-1}\\
&\equiv&g(x)
\end{eqnarray*}
and
\begin{eqnarray*}
g'(x)&=&-2\chp(x)\thp(x)^{p-1}-\shp(x)\cp(x)\thp(x)^{p-2}\\
     & &+\chp(x)\snp(x)^{p-1}\cp(x)^{2-p}+\chp(x)\snp(x)^{p-1}\\
     & &+(p-1)\shp(x)\cp(x)\snp(x)^{p-2}\\
     &\geq&2\chp(x)(\snp(x)^{p-1}-\thp(x)^{p-1})\\
     & &+\shp(x)\cp(x)(\snp(x)^{p-2}-\thp(x)^{p-2})\\
     &\geq&0,
\end{eqnarray*}
where the last inequality follows from $\snp(x)>\thp(x)$. Now it is easy to see that $f(x)>f(0)=0$ which implies the inequality (\ref{2-4.1-l}).
\ep

\begin{conjecture}
For $p\in(2,\infty)$ and $x\in(0,\pi_p/2)$,
\beq\label{2-4.1-l'}
\dfrac{\sinh_p(x)}{x}<\dfrac{p+1}{p+\cos_p(x)}.
\eeq
\end{conjecture}

\bthm
For $p\in[2,\infty)$ and $x\in(0,\pip/2]$,
$$
\dfrac{\snp(x)}{x}>\dfrac{p-1+\cp(x)}{p}\geq\dfrac{1+\cp(x)}{2}.
$$
\ethm

\bp
The second inequality is clear. For the first inequality,
put $f(x)\equiv p\snp(x)-x\cp(x)-(p-1)x$. After some elementary computations, we get
$$
f'(x)=(p-1)\cp(x)+x\cp(x)\tp(x)^{p-1}-(p-1),
$$
and
$$
f''(x)=\cp(x)\tp(x)^{p-2}g(x),
$$
where $g(x)=(p-2)x\tp(x)^p-(p-2)\tp(x)+(p-1)x.$ We have to prove $g(x)>0$ which follows from
$$
g'(x)=p(p-2)x\tp(x)^{p-1}(1+\tp(x)^p)+1>0.
$$
\ep

\blem
For $p>1$, \\
(1) The functions $f_1(x)\equiv\snp(x)/x$  is strictly decreasing from $(0,\pi_p/2)$ onto $(2/\pi_p,1)$. In particular, for $x\in(0,1)$,
\beq
\dfrac{x}{\asp(x)}<\dfrac{\snp(x)}{x}<\dfrac{2x/\pi_p}{\asp(2x/\pi_p)}.
\eeq
(2) The function $f_2(x)\equiv\tp(x)/x$  is strictly increasing from $(0,\pi_p/2)$ onto $(1,\infty)$. In particular, for $x\in(0,k)$,
\beq
\dfrac{x}{\atp(x)}<\dfrac{\tp(x)}{x}<\dfrac{ax}{\atp(ax)},
\eeq
where $0<k<\pi_p/2$ and $a=\tp(k)/k$.\\
(3) The function $f_3(x)\equiv\shp(x)/x$  is strictly increasing from $(0,\infty)$ onto $(1,\infty)$. In particular, for $x\in(0,k)$,
\beq
\dfrac{x}{\ashp(x)}<\dfrac{\shp(x)}{x}<\dfrac{bx}{\ashp(bx)},
\eeq
where $k>0$ and $b=\shp(k)/k$.\\
(4) The function $f_4(x)\equiv\thp(x)/x$  is strictly decreasing from $(0,\infty)$ onto $(0,1)$. In particular, for $x\in(0,k)$,
\beq
\dfrac{x}{\athp(x)}<\dfrac{\thp(x)}{x}<\dfrac{cx}{\athp(cx)},
\eeq
where $k>0$ and $c=\thp(k)/k$.
\elem

\bp
Since the proofs of part (1) to part (4) are similar to each other, we only prove the part (2) here. Since $\tp'(x)=1+\tp(x)^p$ is strictly increasing, the monotone form of l'H\^opital's Rule gives that the function $f_2$ is strictly increasing. Hence we have
$$1<\dfrac{\tp(x)}{x}<\dfrac{\tp(k)}{k}=a,$$
and this is equivalent to
$$\atp(x)<x<\atp(ax).$$
By the monotonicity of $f_2$,
$$\dfrac{x}{\atp(x)}=\dfrac{\tp(\atp(x))}{\atp(x)}<\dfrac{\tp(x)}{x}<\dfrac{\tp(\atp(ax))}{\atp(ax)}=\dfrac{ax}{\atp(ax)}$$
\ep

\bthm
Let $p>1$ and $x>0$. Then \\
(1) $f_1(t)\equiv\cp\left(x/t\right)^t$ is strictly increasing and logarithmic concave in $(2x/\pi_p,\infty)$; \\
(2) $f_2(t)\equiv\snp\left(x/t\right)^t$ is strictly decreasing and logarithmic concave in $(2x/\pi_p,\infty)$; \\
(3) $f_3(t)\equiv\shp\left(x/t\right)^t$ is strictly decreasing and logarithmic concave in $(0,\infty)$;\\
(4) $f_4(t)\equiv\chp\left(x/t\right)^t$ is strictly decreasing and logarithmic convex in $(0,\infty)$.
\ethm

\bp
For part (1), simple computations give
$$\dfrac{d}{dt}\log f_1(t)=\log\cp(s)+s\tp(s)^{p-1}\equiv g_1(s),\quad s=\dfrac{x}{t},$$
and $$g_1'(s)=(p-1)s\tp(s)^{p-2}(1+\tp(s)^p)>0,$$
which implies that $f_1$ is logarithmic concave. For the monotonicity of $f_1$, we write
$h(s)\equiv h_1(s)/h_2(s)$, $h_1(s)\equiv-\log\cp(s)$ and $h_2(s)\equiv s\tp(s)^{p-1}$ with $h_1(0)=h_2(0)=0$, and
$$
\dfrac{h_1'(s)}{h_2'(s)}=\dfrac{1}{1+(p-1)l(s)},\quad l(s)=\dfrac{s(1+\tp(s)^p)}{\tp(s)}=\dfrac{l_1(s)}{l_2(s)}.
$$
By differentiation, we have
$$
\dfrac{l_1'(s)}{l_2'(s)}=1+ps\tp(s)^{p-1}
$$
which is strictly increasing. Hence $h(s)$ is strictly decreasing by the l'H\^opital Monotone Rule, and
$h(s)<h(0)=1/p$ which is equivalent to $s\tp(s)^{p-1}>-p\log\cp(s)$. Now it is easy to see that $g_1(s)>(1-p)\log\cp(s)>0$ which implies that $f_1$ is strictly increasing.

For part (2), it is easy to see that $-t\log(1/\snp(x/t))$ is strictly decreasing in $t$.
Simple computations give
$$\dfrac{d}{dt}t\log\snp(x/t)=-\left(\log\dfrac{1}{\snp(s)}+\dfrac{s}{\tp(s)}\right),\quad s=\dfrac{x}{t},$$
which is strictly increasing in $s$ and hence strictly decreasing in $t$.

For part (3), by differentiations we have
$$\dfrac{d}{dt}\log f_3(t)=\log\shp(s)-\dfrac{s}{\thp(s)}\equiv g_3(s),\quad s=\dfrac{x}{t},$$
and
$$g_3'(s)=\dfrac{s(1-\thp(s)^p)}{\thp(s)^2}>0,$$
which implies that $g_3(s)$ is strictly increasing in $s$ and hence decreasing in $t$. It follows that $\log f_3(t)$ is concave. Since $g_3(s)\leq g_3(\infty)=0$, $f_3$ is strictly decreasing.

For part (4), by differentiations we have
$$\dfrac{d}{dt}\log f_4(t)=\log\chp(s)-s \thp(s)^{p-1}\equiv g_4(s),\quad s=\dfrac{x}{t},$$
and
$$g_4'(s)=-(p-1)s\thp(s)^{p-2}(1-\thp(s)^p)<0,$$
which implies that $g(s)$ is strictly decreasing in $s$ and hence increasing in $t$. It follows that $\log f_4(t)$ is convex. Since $g_4(s)<g_4(0)=0$, $f_4$ is strictly decreasing.
\ep

\begin{openproblem}{\rm
Recently, S. Takeuchi \cite{ta} has introduced functions depending on two parameters $p$ and $q$ that reduce to the functions studied in the present paper when $p=q$.  In \cite{bv1} the authors have continued the study of this family of  generalized  functions, and have suggested that many properties of classical functions also have a counterpart in this more general setting.   It would be natural to generalize the properties of classical trigonometric and hyperbolic functions cited in the survey \cite {avz} to the $(p,q)$-functions of Takeuchi.}
\end{openproblem}

\bigskip

\subsection*{Acknowledgments}
The research of Matti Vuorinen was supported by the Academy of Finland, Project 2600066611.
Xiaohui Zhang is indebted to the Finnish National Graduate School of Mathematics and
its Applications for financial support.


\end{document}